\documentclass{article}
\usepackage[utf8]{inputenc}
\usepackage[left=4cm,right=4cm]{geometry}
\usepackage{cite}
\usepackage{amsmath,amssymb,amsfonts}
\usepackage{graphicx}
\usepackage{textcomp}
\usepackage{upgreek}
\usepackage{bm}
\usepackage{algcompatible}
\usepackage{algorithm}
\usepackage{enumerate}
\usepackage{multirow}
\usepackage{booktabs}

\usepackage{setspace} 
\usepackage{lipsum}
\usepackage[marginal]{footmisc}

\title{Theoretical Study and Comparison of SPSA and RDSA Algorithms with Different Perturbations
}
\author{Yiwen Chen}% <-this % stops a space

\date{}

\begin{document}
\bibliographystyle{plain}
\renewcommand{\baselinestretch}{1.5} 
\normalsize
\maketitle
\newcommand\blfootnote[1]{%
\begingroup
\renewcommand\thefootnote{}\footnote{#1}%
\addtocounter{footnote}{-1}%
\endgroup
}

\blfootnote{1. A compressed version of this paper appears in the 2021 55rd Annual Conference on Information Sciences and Systems (CISS), Mar 2021.\\
2. Contact Information: ychen385@jh.edu / yyiwen1207@gmail.com} 

\begin{abstract}
 Stochastic approximation (SA) algorithms are widely used in system optimization problems when only noisy measurements of the system are available. This paper studies two types of SA algorithms in a multivariate Kiefer-Wolfowitz setting: random-direction SA (RDSA) and simultaneous-perturbation SA (SPSA), and then describes the bias term, convergence, and asymptotic normality of RDSA algorithms. The gradient estimations in RDSA and SPSA have different forms and, consequently, use different types of random perturbations. This paper looks at various valid distributions for perturbations in RDSA and SPSA and then compares the two algorithms using mean-square errors computed from asymptotic distribution. From both a theoretical and numerical point of view, we find that SPSA generally outperforms RDSA.
\end{abstract}

\section{Introduction}

Stochastic approximation is a well-known recursive procedure for finding roots of equations in the presence of noisy measurements. Consider a smooth $p$-dimensional loss function $L:\,\mathbb{R}^p\rightarrow\mathbb{R}$, with gradient $\bm{g}:\mathbb{R}^p \rightarrow\mathbb{R}^p$. Assume that $L$ has a unique minimum $\bm{\uptheta}^*\in\mathbb{R}^p$; i.e. $L(\bm{\uptheta}^*)\leq L(\bm{\uptheta})$ for all $\bm{\uptheta}^*\in\mathbb{R}^p$, and $\bm{g}(\bm{\uptheta}^*)=\bm{0}$. In many cases where only noisy measurements of the gradient are available, the Robbins-Monro stochastic approximation (SA) algorithm is widely used with the form of:
\begin{align}
    \bm{\hat{\uptheta}}_{k+1}=\bm{\hat{\uptheta}}_{k}-a_k\bm{Y}_k(\bm{\hat{\uptheta}}_k).\; k=0,1,2,...
    \label{eq:it}
\end{align}
where $\bm{Y}_k(\bm{\hat{\uptheta}}_k)=\partial{Q}/\partial{\bm{\uptheta}}$ represents a direct noisy measurement of the true gradient $\bm{g}_k(\bm{\hat{\uptheta}}_k)$ for $Q$ in the representation $L(\bm{\uptheta})=E[Q(\bm{\uptheta,V})]$ with $\bm{V}$ corresponding to the randomness in the problem. $a_k>0$ is the step size, determining the convergence of the algorithm.

However, in most practical cases, such direct measurement is difficult to obtain and only noisy measurements of the loss function are available. Then it comes to the use of gradient-free algorithms. There are several common algorithms such as the finite-difference (FD) SA, random-direction SA (RDSA) and the simultaneous-perturbation SA (SPSA). All of them use the measurement of the loss function to estimate the gradient. The FDSA algorithm has also been discussed comprehensively in \cite{spall2005introduction}. It has an $i$th component in the form of:
\begin{align}
    \hat{\bm{g}}_{ki}(\hat{\bm{\uptheta}}_k)=\dfrac{y(\hat{\bm{\uptheta}}_k+c_k\bm{e}_i)-y(\hat{\bm{\uptheta}}_k-c_k\bm{e}_i)}{2c_k},
\label{eq:2}  
\end{align}
where $\bm{e}_i$ denotes the unit vector along the $i$th axis and $y$ is the noisy measurement of the loss value and $c_k>0$ defines the difference magnitude. The pair ${\{a_k,c_k}\}$ represents the gain sequences. Nevertheless, the FDSA algorithm shows inefficiency in solving problems with high dimension due to its  use of $2p$ measurements of $L$ per iteration. The other two RDSA and SPSA approaches, alleviate this problem by requiring only two system simulations regardless of the parameter dimension.

Several papers have discussed the three SA algorithms. Spall has discussed details of FDSA and SPSA algorithms in \cite{717126} and \cite{spall2005introduction}; Kushner discussed the RDSA algorithm in \cite{652378}; Blakney and Zhu compared FDSA and SPSA algorithms in \cite{8693046} and Chin compared the three algorithms using mean-square errors in \cite{558808}. Because the FDSA and SPSA algorithms have been well analyzed, this paper will mainly focus on the RDSA algorithm. Some papers investigated the theoretical foundation for RDSA algorithms such as \cite{book}, \cite{1605408} and \cite{558808}; however, some of them gave the illustration only briefly. This paper is going to arrange and organize those different versions of theoretical proofs to show the convergence and asymptotic normality. After that, we will conduct the comparison between RDSA and SPSA using Mean-Square Error and present certain conditions under which SPSA outperforms RDSA with specific perturbations.

The remainder of this paper is organized as follows. Section 2 introduces the general formulation of the SPSA and RDSA algorithms. Section 3 covers the bias term, convergence condition, and asymptotic normality of RDSA algorithm and leads to the comparison of relative accuracy of SPSA and RDSA. Section 4 presents several numerical studies to assist the theoretical conclusion in Section 3. Especially, this paper uses the skewed-quartic function as an example to show how the theory works on the performance of distinct distributed perturbation in SPSA and RDSA.

\section{Background}

This section briefly discusses the SPSA algorithm and the RDSA algorithm with their estimates for the gradient, $\hat{\bm{g}}(\hat{\bm{\uptheta}}_k)$.

SPSA (`Simultaneous Perturbation Stochastic Algorithm') has been well discussed in \cite[Chap. 7]{spall2005introduction} and \cite{Spall1998ANOO}. The estimate for $\bm{g}(\hat{\bm{\uptheta}}_k)$ is in the form of:
\begin{align*}
    \hat{\bm{g}}_k(\hat{\bm{\uptheta}}_k)=\dfrac{y(\hat{\bm{\uptheta}}_k+c_k\bm{\Delta}_k)- y(\hat{\bm{\uptheta}}_k-c_k\bm{\Delta}_k)}{2c_k\bm{\Delta}_k},
\end{align*}
where $y(\hat{\bm{\uptheta}}_k\pm c_k\bm{\Delta}_k)=L(\hat{\bm{\uptheta}}_k\pm c_k\bm{\Delta}_k)+\upepsilon^{(\pm)}_k$ with $\upepsilon^{(\pm)}_k$ representing measurement noise term that satisfies $E[\upepsilon^{(+)}_k-\upepsilon^{(-)}_k]=0$. $\bm{\Delta}_k\in\mathbb{R}^p$ is a vector of $p$ mutually independent mean-zero random variables representing the perturbation, which satisfies certain conditions, and $\bm{\Delta}^{-1}_k$ represents the vector of inverse components. There's no specific type of distribution for $\bm{\Delta}_k$ as long as the perturbation satisfies general conditions mentioned in \cite[Chap. 7]{spall2005introduction} (Principally, mean 0, symmetric, finite moments and certain inverse moments). Usually, $\bm{\Delta}_k$ can be chosen as Bernoulli, U-shape, and other distributions. The range of choices is restricted due to the main constraint that the inverse moment $E(|1/\Delta_{km}|)$ should be finite. Many common distributions, such as Gaussian distribution and Uniform distribution, cannot be applied to the SPSA algorithm since they contain too much probability mass near 0.

For RDSA (`Random Direction Stochastic Algorithm'), the estimate for $\bm{g}(\hat{\bm{\uptheta}}_k)$ is in the form of:
\begin{align*}
  \hat{\bm{g}}_k(\hat{\bm{\uptheta}}_k)=\dfrac{y(\hat{\bm{\uptheta}}_k+c_k\bm{\pi}_k)- y(\hat{\bm{\uptheta}}_k-c_k\bm{\pi}_k)}{2c_k}\bm{\pi}_k,
\end{align*}
where $\bm{\pi}_k\in\mathbb{R}^p$ is a vector of Monte-Carlo generated random variables satisfying certain regularity conditions and normalized so that $||\bm{\pi}_k||^2=p$ . Unlike SPSA, the requirement for finite inverse moment is released for $\bm{\pi}_k$ in RDSA, so the range of choices of distribution for $\bm{\pi}_k$ is different from that for $\bm{\Delta}_k$ in SPSA. Several valid choices have been discussed in \cite{1605408} such as axis distribution, Gaussian distribution, and Uniformly Spherical distribution where the perturbation distributes uniformly on a $p$-dimensional sphere.

The main difference between the two algorithms is that in RDSA, both adjustments to $\hat{\bm{\uptheta}}_k$ and $\hat{\bm{g}}_k(\hat{\bm{\uptheta}}_k)$ are in the same direction $\bm{\pi}_k$ while the SPSA employs two different directions, $\bm{\Delta}_k$ and $\bm{\upzeta}_k$ with component $\upzeta_{ki}=1/{\Delta}_{ki}$.

There are two popular choices for $\bm{\pi}_k$: independent Gaussian distribution $N(0,1)$ and Uniformly Spherical distribution with radius $\sqrt{p}$. Notice that in Uniformly Spherical distribution, since there is the normalization, $||\bm{\pi}_k||^2=p$, variables are not independent, but uncorrelated.

\section{Methodology}

This section presents several results that form the theoretical basis for the RDSA algorithm, similar to those results of SPSA, which have been completely proved in \cite{spall1992multivariate}. The following sections consider the bias term in $\bm{\hat{g}}_k(\bm{\hat{\uptheta}}_k)$, the strong convergence and asymptotic normality of $\hat{\bm{\uptheta}}_k$ in the RDSA algorithm. Then we compare the performance of RDSA with SPSA using mean-square error computed with the asymptotic distribution.

\subsection{The bias in $\hat{\bm{g}}_k$}

The bias term in $\hat{\bm{g}}_k$ is in the form of:
\begin{align}
   \bm{b}_k(\hat{\bm{\uptheta}}_k)=E[\hat{\bm{g}}_k(\hat{\bm{\uptheta}}_k)|\hat{\theta}_k]-\bm{g}_k(\bm{\hat{\uptheta}}_k),  
\label{eq:6}
\end{align}
where $\upepsilon^{(\pm)}_k=y(\bm{\uptheta}_k\pm c_k\bm{\pi}_k)-L(\bm{\uptheta}_k\pm c_k\bm{\pi}_k)$ is the noise on the estimate of the loss function $L$, with mean 0 and variance $\upsigma^2$. Here we henceforth assume that the noise is independent with $\hat{\bm{\uptheta}}_k$ and $\bm{\pi}_k$.

Calculate $E[\hat{\bm{g}}_k(\hat{\bm{\uptheta}}_k)]$ and $\bm{g}(\hat{\bm{\uptheta}}_k)$ respectively:
\begin{align*}
    \hat{\bm{g}}_k(\hat{\bm{\uptheta}}_k)=&\dfrac{y(\hat{\bm{\uptheta}}_k+c_k\bm{\pi}_k)-y(\hat{\bm{\uptheta}}_k-c_k\bm{\pi}_k)}{2c_k}\;\bm{\pi}_k\\
    &=\dfrac{L(\hat{\bm{\uptheta}}_k+c_k\bm{\pi}_k)+\upepsilon^{(+)}_k-L(\hat{\bm{\uptheta}}_k-c_k\bm{\pi}_k)-\upepsilon^{(-)}_k}{2c_k}\;\bm{\pi}_k\\
    &=\dfrac{L(\hat{\bm{\uptheta}}_k+c_k\bm{\pi}_k)-L(\hat{\bm{\uptheta}}_k-c_k\bm{\pi}_k)}{2c_k}\;\bm{\pi}_k+\dfrac{\upepsilon^{(+)}_k-\upepsilon^{(-)}_k}{2c_k}\bm{\pi}_k\\
    &=\dfrac{2c_k\bm{g}(\hat{\bm{\uptheta}}_k)^T\bm{\pi}_k+\dfrac{1}{6}c^3_kL^{'''}(\bar{\bm{\uptheta}}^{(+)}_k)[\bm{\pi}_k\otimes\bm{\pi}_k\otimes\bm{\pi}_k]}{2c_k}\;\bm{\pi}_k\\
    &+\dfrac{\dfrac{1}{6}c^3_kL^{'''}(\bar{\bm{\uptheta}}^{(-)}_k)[\bm{\pi}_k\otimes\bm{\pi}_k\otimes\bm{\pi}_k]}{2c_k}\;\bm{\pi}_k\\
    &+O(c^5_k)+\dfrac{\upepsilon^{(+)}_k-\upepsilon^{(-)}_k}{2c_k}\bm{\pi}_k
\end{align*}

where $\bar{\bm{\uptheta}}^{(\pm)}_k$ denotes points on the line segments between $\hat{\bm{\uptheta}}_k$ and $\hat{\bm{\uptheta}}_k\pm c_k\bm{\pi}_k$.

Then look at the $m$th term of $\hat{\bm{g}}_k(\hat{\bm{\uptheta}}_k)$:
\begin{align*}
    \hat{g}_{km}(\hat{\bm{\uptheta}}_k)=&\left(\sum_{i=1}^{p}g_i(\hat{\bm{\uptheta}}_k){\pi}_{ki}+\dfrac{1}{12}c^2_k[L^{'''}(\bar{\bm{\uptheta}}^{(+)}_k)+L^{'''}(\bar{\bm{\uptheta}}^{(-)}_k)][\bm{\pi}_k\otimes\bm{\pi}_k\otimes\bm{\pi}_k]\right)\,\pi_{km}\\
    &+\dfrac{\upepsilon^{(+)}_{km}-\upepsilon^{(-)}_{km}}{2c_k}\pi_{km}\\
    &=g_m(\hat{\uptheta}_k)\pi^2_{km}+\sum_{i=1,i\neq m}^pg_i(\hat{\uptheta}_k)\pi_{ki}{\pi}_{km}\\
    &+\dfrac{1}{12}c^2_k\left\{[L^{'''}(\bar{\bm{\uptheta}}^{(+)}_k)+L^{'''}(\bar{\bm{\uptheta}}^{(-)}_k)][\bm{\pi}_k\otimes\bm{\pi}_k\otimes\bm{\pi}_k])\;\pi_{km}\right\}+\dfrac{\upepsilon^{(+)}_{k}-\upepsilon^{(-)}_{k}}{2c_k}\pi_{km}
\end{align*}

Based on the basic assumption of independence between $\upepsilon^{(\pm)}_k$ and $\bm{\pi}_k$ and the uncorrelation between $\pi_{ki}$ and $\pi_{kj}$, $i\neq j$, we have: $$E\left[\dfrac{\upepsilon^{(+)}_k-\upepsilon^{(-)}_k}{2c_k}\bm{\pi}_k\right]=0;\;E(\pi^2_{km})=1;\; E(\pi_{ki}\pi_{km})=0.$$ 
Then,
\begin{small}
\begin{align*}
    E[\hat{g}_{km}(\hat{\bm{\uptheta}}_k)]&=g_m(\hat{\bm{\uptheta}}_k)E(\pi^2_{km})+\sum_{i=1,i\neq m}^pg_i(\hat{\bm{\uptheta}}_k)E(\pi_{ki}\pi_{km})\\
    &+\dfrac{c^2_k}{12}E\left\{[L^{'''}(\bar{\bm{\uptheta}}^{(+)}_k)+L^{'''}(\bar{\bm{\uptheta}}^{(-)}_k)][\bm{\pi}_k\otimes\bm{\pi}_k\otimes\bm{\pi}_k]\;\pi_{km}\right\}\\
    &=g_m(\hat{\bm{\uptheta}}_k)+\dfrac{c^2_k}{12}E\left\{[L^{'''}(\bar{\bm{\uptheta}}^{(+)}_k)+L^{'''}(\bar{\bm{\uptheta}}^{(-)}_k)][\bm{\pi}_k\otimes\bm{\pi}_k\otimes\bm{\pi}_k]\pi_{km}\right\}\\
    &=g_m(\hat{\bm{\uptheta}}_k)+\dfrac{1}{6}c^2_kE\left\{ [L^{'''}(\hat{\bm{\uptheta}}_k)+O(c^2_k)][\bm{\pi}_k\otimes\bm{\pi}_k\otimes\bm{\pi}_k]\;\pi_{km}\right\}\\
    &=g_m(\hat{\bm{\uptheta}}_k)+\dfrac{1}{6}c^2_k\left(\sum_{i,j,k}\dfrac{\partial^3L}{\partial\pi_{ki}\partial\pi_{kj}\partial\pi_{kl}}\;\pi_{ki}\pi_{kj}\pi_{kl}\right)\;\pi_{km}+O(c^4_k)\\
    &=g_m(\hat{\bm{\uptheta}}_k)+\dfrac{1}{6}c^2_k\left(\dfrac{\partial^3L(\hat{\bm{\uptheta}}_k)}{\partial{\pi}_{km}^3}\;\pi^4_{km}+\sum_{i,j,l\neq m}\dfrac{\partial^3L}{\partial\pi_{ki}\partial\pi_{kj}\partial\pi_{kl}}\;\pi_{ki}\pi_{kj}\pi_{kl}\pi_{km}\right)+O(c^4_k)
\end{align*}
\end{small}
So, the $m$th bias term is in the form of:
\begin{align*}
\begin{split}
    \bm{b}_{km}(\hat{\bm{\uptheta}}_k)&=E[\hat{g}_{km}(\hat{\bm{\uptheta}}_k)]-g_m(\hat{\bm{\uptheta}}_k)\\
    &=\dfrac{1}{6}c^2_k\left(L^{'''}_{mmm}(\hat{\bm{\uptheta}}_k)E(\pi^4_{km})+3\sum\limits_{i=1,i\neq m}^pL^{'''}_{iim}(\hat{\bm{\uptheta}}_k)E(\pi^2_{ki}\pi^2_{km})\right)+O(c^4_k)
\end{split}
\end{align*}

\subsection{Convergence Conditions}

This subsection presents conditions for convergence of the RDSA iteration: $\bm{\uptheta}^*\rightarrow\bm{\uptheta}^*\; \text{a.s.}$. The following conditions ensure convergence of $\hat{\bm{\uptheta}}_k$ to $\bm{\uptheta}^*$ in the RDSA algorithm:

Consider all $k\geq K$ for some $K<\infty$. Suppose that for each $k$,  there exists some $\upalpha_0,\upalpha_1,\upalpha_2>0$, $\bm{\pi}_{k}$ is a direction vector normalized so that $||\bm{\pi}_k||^2=p$ and symmetrically distributed about 0 with $|\pi_{ki}^4|\leq\upalpha_0\;\text{a.s. and}\; E|\pi_{ki}|^2\leq\upalpha_1. i=1,2,...,p $ For almost all $\hat{\bm{\uptheta}}_k$, suppose that for $\forall\;\bm{\uptheta}$ in an open neighborhood of $\hat{\bm{\uptheta}}_k$, $L^{'''}(\bm{\uptheta})$ exists continuously with individual elements satisfying $|L^{'''}_{ijk}(\bm{\uptheta})|\leq\upalpha_2$. Similar to the conditions for SPSA discussed in \cite[Chap.7]{spall2005introduction}, if the derivatives of $L$ are equicontinuous and bounded, then the RDSA algorithm will converge almost surely under the following conditions:

\begin{itemize}
    \item A1: $a_k, c_k >0,\;\forall k: a_k\rightarrow0, c_k\rightarrow0\;{\text{as}}\; k\rightarrow\infty,\;\sum\limits^{\infty}_{k=0}a_k=\infty,\;\sum\limits^{\infty}_{k=0}\left(\dfrac{a_k}{c_k}\right)^2<\infty.$
    \item A2: $\mathop{{\text{sup}}}\limits_k||\hat{\bm{\uptheta}}_k||<\infty\; \text{a.s.}$
    \item A3: $\bm{\uptheta}^*$ is an asymptotically stable solution of the differential equation: $d\bm{x}(t)/dt=-\bm{g}(\bm{x}).$
    \item A4: Let $D(\bm{\uptheta}^*)=[x_0|\mathop{\text{lim}}\limits_{t\rightarrow\infty}\bm{x}(t|x_0)=\bm{\uptheta}^*]$ where $\bm{x}(t|x_0)$ denotes the solution to the differential equation of $d\bm{x}(t)/dt=-\bm{g}(\bm{x})$ based on initial conditions $\bm{x}_0$. There exists a compact $S\subseteq D(\bm{\uptheta}^*)$ s.t $\Tilde{\bm{\uptheta}}\in S$ infinitely often for almost all sample points.
    \item A5: $\forall k, E|\upepsilon^{(\pm)}_k| ^2\leq\upalpha_0,E(\bm{\pi}_k\bm{\pi}^T_k)=\bm{I}$, and $E|\pi^2_{kl}L(\hat{\bm{\uptheta}}_k\pm c_k{\bm{\pi}}_k)^2|\leq\upalpha_1, \;l=1,2,...,p.$
\end{itemize}

\textbf{Proposition 1A:} Let A1--A5 hold and for some $K<\infty$, each $k\geq K$, ${\{\pi_{ki}\}}$ is Gaussian (0,1) distributed i.i.d., then $\hat{\bm{\uptheta}}_k\rightarrow\bm{\uptheta}^*$ as $k\rightarrow\infty$.

\textbf{Proposition 1B:} Let A1--A5 hold and for some $K<\infty$, each $k\geq K$, ${\{\pi_{ki}\}}$ are uniformly distributed on a $p$-dim sphere, then $\hat{\bm{\uptheta}}_k\rightarrow\bm{\uptheta}^*$ as $k\rightarrow\infty$.

\noindent {\it{Proof:}} This proof applies for both 1A and 1B above. From the above subsection in terms of bias in $\hat{\bm{g}}_k$ and the A1 condition, we know that:
\begin{align}
 ||b_k(\hat{\bm{\uptheta}}_k)||<\infty \qquad \forall k,\;b_k(\hat{\bm{\uptheta}}_k)\rightarrow\bm{0}\;\;\text{a.s.}
\label{eq:11}
\end{align}

According to the inequality in \cite[p. 315]{bookStochastic}
, it holds that:
\begin{align}
%\begin{split}
    &P(\mathop{\text{sup}}\limits_{m\geq k}||\sum\limits_{i=k}^ma_i\bm{e}_i||\geq\upeta)\leq\upeta^{-2}E||\sum\limits_{i=k}^\infty a_i\bm{e}_i||^2=\upeta^{-2}\sum\limits_{i=k}^\infty a^2_iE||\bm{e}_i||^2.
\label{eq:12}
%\end{split}
\end{align}

By the definition of $\bm{e}_k(\hat{\bm{\uptheta}}_k)$ and Condition A1, A5, we would have:$$\mathop{\text{lim}}\limits_{k\rightarrow\infty}P(\mathop{\text{sup}}\limits_{m\geq k}||\sum\limits_{i=k}^ma_i\bm{e}_i||\geq\upeta)=0\;\; \text{for\;any\;} \upeta >0.$$

Then, by \cite[Lemma 2.2.1]{10.2307/2239701}, the above propositions hold.

\subsection{Asymptotic Normality}

This subsection mainly discusses the asymptotic normality for $\hat{\bm{\uptheta}}_k$ in RDSA algorithm, especially, how the Fabian's Theorem in \cite[Theorem 2.2]{10.2307/2239701} applies on RDSA algorithm based on its generalization for SPSA shown in \cite{spall1992multivariate} and \cite{HERNANDEZ2019420}. In order to satisfy the conditions for Fabian's Theorem, we strengthen the Condition A5 to:

A5$^\prime$: For some $\updelta, \upalpha_0, \upalpha_1, \upalpha_2>0$ and $\forall k, E|\upepsilon^{(\pm)}_k|^{2+\delta} \leq\upalpha_0, E(\pi_{kl}L(\hat{\bm{\uptheta}}_k\pm c_k\bm{\pi}_k))^{2+\updelta}\leq\upalpha_1, E(\pi_{kl})^{4+\updelta}\leq\upalpha_2,\;l=1,2,...,p,  E(\bm{\pi}_k\bm{\pi}^T_k)=\bm{I}$. 

Then, the conditions for the asymptotic normality for $\hat{\bm{\uptheta}}_k$ in RDSA can be stated as the following proposition:

\textbf{Proposition 2:} Assume that conditions A1--A4 and A5$^\prime$ hold and $\upbeta>0$. Let $\bm{H}(\cdot)$ denote the Hessian matrix for $L(\bm{\uptheta})$; $\upsigma, \upphi$ be such that $E(\upepsilon_k^{(+)}-\upepsilon_k^{(-)})^2\rightarrow\upsigma^2$ and $E(\pi_{kl}^4)\rightarrow\upphi,\; k\rightarrow\infty$. Suppose the gain sequences are $a_k=a/(k+1)^\upalpha, c_k=c/(k+1)^\upgamma, a>0, c>0, k=0,1,2,3,...$ and set $0<\upalpha\leq1, \upgamma\geq\upalpha/6, \upbeta=\upalpha-2\upgamma$. Let $\bm{P}$ be an orthogonal matrix such that $\bm{PH}(\bm{\uptheta}^*)\bm{P}^T=a^{-1}\text{diag}(\uplambda_1,...,\uplambda_p)$. Then
\begin{align}
   k^{\upbeta/2}(\hat{{\uptheta}}_k-{\uptheta}^*)\stackrel{\text{dist}}{\longrightarrow}N(\bm{\upmu},\bm{PMP}^T),\quad k\rightarrow\infty, \label{eq:notmal}
\end{align}
where $\bm{M}=\dfrac{1}{4c^2}a^2\upsigma^2\text{diag}[1/(2\uplambda_1-\upbeta_{+}),....,1/(2\uplambda_p-\upbeta_{+})]$ with $\upbeta_{+}=\upbeta<2\mathop{\text{min}}\limits_{i}\uplambda_i$ if $\upalpha=1$ and $\upbeta_{+}=0$ if $\upalpha<1$, and 
\begin{align}
    \bm{\upmu}=\left\{
\begin{array}{ccr}
0 & & {\text{if}\;3\upgamma-\upalpha/2>0}\\
(a\bm{H}(\uptheta^*)-\dfrac{1}{2}\upbeta_{+}\bm{I})^{-1}\bm{T} &  & {\text{if}\;3\upgamma-\upalpha/2=0}\\
\end{array}
\right.,
\label{eq:mu}
\end{align}
where the $l$th component of $\bm{T}$ is:
\begin{align*}
    \bm{T}^{\text{RD}}_l=-\dfrac{1}{6}ac^2\left[{\upphi} L^{'''}_{lll}(\bm{\uptheta}^*)+3{\upsilon}\sum\limits_{i\neq l}^pL^{'''}_{iil}(\bm{\uptheta}^*)\right],
\end{align*}
and
\begin{align*}
   {\upphi}=E(\pi^4_{kl}),\;{\upsilon}=E(\pi^2_{kl}\pi^2_{km}),\;l\neq m.    
\end{align*}

\noindent {\it{Proof:}} The result will be shown if conditions (2.2.1), (2.2.2), and (2.2.3) of \cite[Theorem 2.2]{10.2307/2239701} hold. In the notation of \cite[Theorem 2.2]{10.2307/2239701}, we can rewrite the formula (\ref{eq:it}) in the form of:
\begin{align*}
%\begin{split}
    \hat{\bm{\uptheta}}_{k+1}-\bm{\uptheta}^*=(\bm{I}-k^{-\upalpha}\bm{\Gamma}_k)(\hat{\bm{\uptheta}}_k-\bm{\uptheta}^*)+k^{-(\upalpha+\upbeta)/2}\bm{\Phi}_k\bm{V}_k+k^{-\upalpha-\upbeta/2}\bm{T}_k,
%\end{split}
\end{align*}
where 
\begin{align*}
     \bm{\Gamma}_k &=k^{\upalpha}a_k\Tilde{\bm{H}}_k,\;\;\;\bm{\Phi}_k=-a\bm{I},\\ 
     \bm{V}_k & =\dfrac{1}{k^\upgamma}[\hat{\bm{g}}_k(\hat{\bm{\uptheta}}_k)-E(\hat{\bm{g}}_k(\hat{\bm{\uptheta}}_k)],\;\;\bm{T}_k=-ak^{\upbeta/2}\bm{b}_k(\hat{\bm{\uptheta}}_k).
\end{align*}

Next, we are going to verify Fabian's conditions.

As shown in the \cite[p. 233]{cuevas2017cyclic}, let $\bm{\Gamma}_k=k^{\upalpha}a_k\Tilde{\bm{H}}_k$ where the $i$th row of $\Tilde{\bm{H}}_k$ equals to the $i$th row of the Hessian matrix evaluated at $\bm{\uptheta}=(1-\uplambda_i)\hat{\bm{\uptheta}}_k+\uplambda_i\bm{\uptheta}^*$ for some $\uplambda_i\in[0,1]$. Due to the continuity of $\bm{H}(\bm{\uptheta})$ at $\bm{\uptheta}^*$ and $k^{\upalpha}a_k\rightarrow a, k\rightarrow\infty$, $a\bm{H}(\bm{\uptheta}^*)$ would be a positive definite matrix, which $\bm{\Gamma}_k$ would converge to.

Next, consider the convergence of $\bm{T}_k$:
\begin{align*}
    T_{kl}&=-\dfrac{1}{12}ac^2\dfrac{1}{k^{3\upgamma-\upalpha/2}}E\left\{\pi_{kl}[L^{'''}(\bar{\bm{\uptheta}}^{(+)}_k)+L^{'''}(\bar{\bm{\uptheta}}^{(-)}_k)]\bm{\pi}_k\otimes\bm{\pi}_k\otimes\bm{\pi}_k\right\}.
\end{align*}

If $3\upgamma-\upalpha/2>0$, we have ${T}_{kl}\rightarrow0$ a.s. 

If $3\upgamma-\upalpha/2=0$, use the fact that $L^{'''}$ is uniformly bounded near $\bm{\uptheta}^*$, we have:
\begin{small}
$$
    T_{kl}\rightarrow-\dfrac{ac^2}{6}\left\{E(\pi^4_{kl})L^{'''}_{lll}(\bm{\uptheta}^*)+3\sum\limits_{i\neq l}^p L^{'''}_{iil}(\bm{\uptheta}^*)E(\pi^2_{ki}\pi^2_{kl})\right\}\; \text{a.s.}
$$
\end{small}

So, we can show that $\bm{T}_k$ converges a.s. for $3\upgamma-\upalpha/2\geq0.$. Further, obviously, $\bm{\Phi}_k=-a\bm{I}\rightarrow\bm{\Phi}=-a\bm{I}.$

Last, consider $\bm{V}_k$:
\begin{small} 
    \begin{align}
    \begin{split}
    &E(\bm{V}_k\bm{V}^T_k)\\
    &=k^{-2\upgamma}E\bigg\{\bm{\pi}_k(\bm{\pi}_k)^T\left[\dfrac{L(\hat{\bm{\uptheta}}_k+c_k\bm{\pi}_k)-L(\hat{\bm{\uptheta}}_k-c_k\bm{\pi}_k)}{2ck^{-\upgamma}}\right]^2\\
    &+ \bm{\pi}_k(\bm{\pi}_k)^T\left[\dfrac{\upepsilon^{(+)}_k-\upepsilon^{(-)}_k}{2ck^{-\upgamma}}\right]\left[\dfrac{L(\hat{\bm{\uptheta}}_k+c_k\bm{\pi}_k)-L(\hat{\bm{\uptheta}}_k-c_k\bm{\pi}_k)}{2ck^{-\upgamma}}\right]\\
    &+\bm{\pi}_k(\bm{\pi}_k)^T\left[\dfrac{\upepsilon^{(+)}_k-\upepsilon^{(-)}_k}{2ck^{-\upgamma}}\right]^2\bigg\}\\
    &-k^{-2\upgamma}\left[\bm{g}(\hat{\bm{\uptheta}}_k)+\bm{b}_k(\hat{\bm{\uptheta}}_k)\right]\left[\bm{g}(\hat{\bm{\uptheta}}_k)+\bm{b}_k(\hat{\bm{\uptheta}}_k)\right]^T.
    \label{eq:chang}
    \end{split}
\end{align}
\end{small}

Similar to the proof (3.5) in \cite{spall1992multivariate}, for sufficiently large $k$, $L(\hat{\bm{\uptheta}}_k\pm c_k{\bm{\pi}_k})$ is uniformly bounded in $\bm{\pi}_k$. Combined with condition A2 and Holder's Inequality, it is implied that the first and second terms in (\ref{eq:chang}) will converge to 0 a.s. Also, by conditions A1--A5, the fourth term will converge to 0 a.s. as well.

Then consider the third term:
\begin{align*}
    &E\left\{\bm{\pi}_k(\bm{\pi}_k)^T\left[\upepsilon^{(+)}_k-\upepsilon^{(-)}_k\right]^2\right\}\\
    &=\int_{\Omega_{\bm{\pi}}}\bm{\pi}_k(\bm{\pi}_k)^TE\left[\upepsilon^{(+)}_k-\upepsilon^{(-)}_k\right]^2dP_{\bm{\pi}}.
\end{align*}
where $\Omega_{\bm{\pi}}$ is the sample space generating the $\bm{\pi}_k$ and $P$ is the corresponding probability measure. Here we have the fact that $E\left[\upepsilon^{(+)}_k-\upepsilon^{(-)}_k\right]^2\rightarrow\upsigma^2$ a.s. 
and $||\bm{\pi}_k||^2=p$. If we consider ${\{\bm{\pi}_{ij}}\}$ to be independently generated with mean 0 and second moment 1 and then it leads to the result that $E({\pi}_{kl}{\pi}_{km})=E({\pi}_{kl})E({\pi}_{km})=0$. If we consider $\{\pi_k\}$ being generated from a $p$-dimensional sphere with radius $\sqrt{p}$ uniformly, then for $l\neq m,\;{\pi}_{kl},{\pi}_{km}$ are not independent, but uncorrelated. In this case,  we still have $E({\pi}_{kl}{\pi}_{km})=E({\pi}_{kl})E({\pi}_{km})=0$. Moreover, $E({\pi}^2_{kl})=E({\pi}^2_{km})$, such that $E({\pi}^2_{kl})=1$. So, we can have: $E[\bm{\pi}_k(\bm{\pi}_k)^T]=\bm{I},\;\;\forall \;l\neq m$. such that,$$E(\bm{V}_k\bm{V}^T_k)\rightarrow\dfrac{1}{4}c^{-2}\upsigma^2\bm{I}.$$

Now, the conditions (2.2.1) and (2.2.2) in \cite[Theorem 2.2]{10.2307/2239701} have been shown. Next, we are going to show the condition (2.2.3) also holds. By the definition of $\bm{V}_k$, Markov Inequality, and the triangle inequality and the proof shown in (3.8) in \cite{spall1992multivariate}, for any $0<\upsilon<\updelta/2$, it holds that:
\begin{align*}
   ||\bm{V}_k||^{2(1+\upsilon)}\leq&2^{2(1+\upsilon)}k^{-2\upgamma(1+\upsilon)}[||\hat{\bm{g}}_k(\hat{\bm{\uptheta}}_k)||^{2(1+\upsilon)}\\
   &+||\bm{b}_k(\hat{\bm{\uptheta}}_k)||^{2(1+\upsilon)}+||\bm{g}_k(\hat{\bm{\uptheta}}_k)||^{2(1+\upsilon)}].    
\end{align*}

Since $\bm{g}(\hat{\bm{\uptheta}}_k)$ and $\bm{b}_k(\hat{\bm{\uptheta}}_k)$ are uniformly bounded and $L(\hat{\bm{\uptheta}}_k\pm c_k{\bm{\pi}}_k)$ is uniformly bounded for $\forall k>K$, we have
\begin{align*}
    &E[k^{-2\upgamma(1+\upsilon)}||\bm{b}_k(\hat{\bm{\uptheta}}_k)||^{2(1+\upsilon)}]\rightarrow\bm{0},\\
    &E[k^{-2\upgamma(1+\upsilon)}||\bm{g}_k(\hat{\bm{\uptheta}}_k)||^{2(1+\upsilon)}]\rightarrow\bm{0}.
\end{align*}
as $k\rightarrow\infty$. Moreover, invoking A5$^\prime$, Holder's Inequality, similar to the case in \cite{spall1992multivariate}, it implies that $E[||\hat{g}_k(\hat{\bm{\uptheta}}_k)||^{2(1+\upsilon)}]=O(k^{2\upgamma(1+\upsilon)})$. As a result, $E||\bm{V}_k||^{2(1+\upsilon)}=O(1)$, which shows that: $$\mathop{\rm{lim}}\limits_{k\rightarrow\infty}E\left(I_{||\bm{V}_k||^2\geq rk^{\upalpha}}||\bm{V}_k||^2\right)=0,\qquad \forall r>0,$$
where $I_{\{.\}}$ is the indicator function.

Now all the required conditions for Fabian's Theorem in \cite[Theorem 2.2]{10.2307/2239701} have been verified. It is stated that $\hat{\bm{\uptheta}}_k$ in RDSA algorithm is asymptotically normal.

\textbf{Corollary 2A:} Assume conditions in Proposition 2 hold and perturbation $\{{\bm{\pi}_k\}}$ is Gaussian distributed $N(0,1)$ independently. Then the $l$th component of $\bm{T}$ is:
$$\bm{T}^{\text{RD}}_l=-\dfrac{1}{6}ac^2\left[3 L^{'''}_{lll}(\bm{\uptheta}^*)+3\sum\limits_{m\neq l}^pL^{'''}_{llm}(\bm{\uptheta}^*)\right].$$

\textbf{Corollary 2B:} Assume conditions in Proposition 2 hold and perturbation $\{{\bm{\pi}_k}\}$ is uniformly distributed on a $p$-dim sphere with radius $\sqrt{p}$. Then the $l$th element of $\bm{T}$ is:
$$\bm{T}^{\text{RD}}_l=-\dfrac{1}{6}ac^2\left[\dfrac{3p}{p+2} L^{'''}_{lll}(\bm{\uptheta}^*)+\dfrac{3p}{p+2}\sum\limits_{m\neq l}^pL^{'''}_{llm}(\bm{\uptheta}^*)\right].$$

\subsection{Relative Accuracy between SPSA and RDSA}
\subsubsection{Introduction}

This subsection discusses the MSE for SPSA and RDSA algorithm using the same gain sequence. We mainly focus on several widely used distributions for perturbations in SPSA and RDSA. We choose Bernoulli, and U-shape, especially polynomial with high order for $\bm{\Delta}_k$. (Here we choose $\Delta_{ki}\stackrel{\text{i.i.d}}\sim x^{10},\,x\in[-1.17,1.17]$ to make it a valid probability density function) and choose Gaussian and uniform-spherical for $\bm{\pi}_k$ in RDSA.

Notice that $\hat{{\bm{\uptheta}}}_k$ in the SPSA and RDSA algorithms have the asymptotic normality in, respectively, \cite{spall1992multivariate} or formula (\ref{eq:notmal}) with different values for the parameters. Here we give several notations to make the parameters more concise: 
\begin{align*}
\begin{split}
    &\upphi=E[({\pi}^4_{ki})]\;\;;\;\,\qquad {\upsilon}=E({\pi}^2_{ki}{\pi}^2_{km}),\;i\neq m\\
    &{\uprho}^2=E\left[({\Delta}_{ki})^{-2}\right];\;\;{\xi}^2=E\left[({\Delta}_{ki})^2\right].
    \end{split}
\end{align*}

Then we have:
\begin{align*}
    \bm{T}^{\text{RD}}_l &=-\dfrac{1}{6}ac^2\left[{\upphi} L^{'''}_{lll}(\bm{\uptheta}^*)+3{\upsilon}\sum\limits_{m\neq l}^pL^{'''}_{llm}(\bm{\uptheta}^*)\right];\\
    \bm{M}^{\text{RD}} &=\dfrac{1}{4c^2}a^2\upsigma^2\text{diag}\left(\dfrac{1}{2\uplambda_1-\beta_{+}},....,\dfrac{1}{2\uplambda_p-\beta_{+}}\right);\\
    \bm{T}^{\text{SP}}_l &=-\dfrac{1}{6}ac^2{\xi}^2\left[L^{'''}_{lll}(\bm{\uptheta}^*)+3\sum\limits_{m\neq l}^pL^{'''}_{llm}(\bm{\uptheta}^*)\right];\\
    \bm{M}^{\text{SP}} &=\dfrac{1}{4c^2}a^2\upsigma^2{\uprho}^2\text{diag}\left(\dfrac{1}{2\uplambda_1-\upbeta_{+}},....,\dfrac{1}{2\uplambda_p-\beta_{+}}\right).
\end{align*}

Notice that tr$(\bm{M}_{\text{SP}})={\uprho}^2$tr$(\bm{M}_{\text{RD}}).$ Then, we can compute the asymptotic mean-square error by the definition:
\begin{align*}
   {\text{MSE}}=\bm{\upmu}^T\bm{\upmu}+{\text{tr}}(\bm{PMP}^T)=\bm{\upmu}^T\bm{\upmu}+{\text{tr}}(\bm{M}).
\end{align*}

Asymptotically, the ratio of the MSE values of the two algorithms are in the form of:
 \begin{small}
\begin{align}
\begin{split}
    \dfrac{\text{MSE}_\text{RD}}{\text{MSE}_\text{SP}}&\rightarrow\dfrac{\bm{\mu}^T_{\text{RD}}\bm{\mu}_{\text{RD}}+{\text{tr}}(\bm{PM}_{\text{RD}}\bm{P}^T)}{\bm{\mu}^T_{\text{SP}}\bm{\mu}_{\text{SP}}+{\text{tr}}(\bm{PM}_{\text{SP}}\bm{P}^T)}\\
    &={\left[{\upphi} \bm{u}_1+{\upsilon}\bm{u}_2\right]^T\bm{S}\left[{\upphi} \bm{u}_1+{\upsilon}\bm{u}_2\right]+{D}}{\left[\upxi^2 \bm{u}_1+\upxi^2\bm{u}_2\right]^T\bm{S}\left[\upxi^2 \bm{u}_1+\upxi^2\bm{u}_2\right]+{\uprho}^2{D}},
\label{eq:ss}
\end{split}
\end{align}
\end{small}
where 
\begin{align*}
    &\bm{u}_1=ac^{2}[L^{'''}_{111}(\bm{\uptheta}^*),...,L^{'''}_{ppp}(\bm{\uptheta}^*)]^T/6;\\
    &\bm{u}_2=ac^{2}[3\sum\limits_{j\neq 1}L^{'''}_{jj1}(\bm{\uptheta}^*),..., 3\sum\limits_{j\neq p}L^{'''}_{jjp}(\bm{\uptheta}^*)]^T/6;\\
    &\bm{S}=\left(a\bm{H}(\bm{\uptheta}^*)-\dfrac{\upbeta_{+}}{2}\bm{I}\right)^{-2};\\
    &D={\text{tr}}(\bm{M}_{\text{RD}})=\dfrac{a^2\upsigma^2}{4c^2}\sum\limits_{i=1}^p\dfrac{1}{2\uplambda_i-\upbeta_{+}}.
\end{align*}

Denote
\begin{align*}
    \begin{split}
        &{Q}_1=(\bm{u}_1+\bm{u}_2)^T\bm{S}(\bm{u}_1+\bm{u}_2),\\ &{Q}_2=(3\bm{u}_1+\bm{u}_2)^T\bm{S}(3\bm{u}_1+\bm{u}_2).
    \end{split}
\end{align*}

The parameters ${\{\upphi,\upsilon,\upxi^2,\uprho^2}\}$ and the MSE values for several common used perturbations in SPSA and RDSA are listed as follows:

\begin{table}[ht!]
 \renewcommand\arraystretch{1.5}
 \centering
 \caption{Key parameters for perturbations in RDSA and SPSA algorithms with different distributions}
  \begin{tabular}{cccc}
  \toprule  %添加表格头部粗线
  \small{Perturbation in RDSA} &${\upphi}$&${\upsilon}$\\
  \midrule  %添加表格中横线
   ${\pi}_{ki}\stackrel{\text{i.i.d}}\sim N(0,1)$&3&1\\
   $\bm{\pi}_k\sim$ Uniformly Spherical &$\dfrac{3p}{p+2}$&$\dfrac{p}{p+2}$\\
  \bottomrule %添加表格底部粗线
  \end{tabular}
  \label{table1}
\end{table}

\begin{table}[ht!]
 \renewcommand\arraystretch{1.5}
 \centering

  \begin{tabular}{ccc}
  \toprule  %添加表格头部粗线
  \small{Perturbation in SPSA} &${\xi}^2$&${\rho}^2$\\
  \midrule  %添加表格中横线
  ${\Delta}_{ki}\stackrel{\text{i.i.d}}\sim$ Bernoulli ($\pm1$)&1&1\\
  ${\Delta}_{ki}\stackrel{\text{i.i.d}}\sim x^{10}$,\, $x\in[-1.17,1.17]$
  & 1.15 &0.90\\
  \bottomrule %添加表格底部粗线
  \end{tabular}
  \label{table2}
\end{table}

\begin{table}[ht!]
 \renewcommand\arraystretch{1.5}

 \centering
 \caption{MSE in SPSA and RDSA with different perturbations}
  \begin{tabular}{ccc}
  \toprule  %添加表格头部粗线
   Algorithm & Distribution & MSE\\
  \midrule  %添加表格中横线
  SPSA & Bernoulli&${Q}_1+D$\\
  SPSA & U-shape ($x^{10}$) &$1.33{Q}_1 + 0.9D$\\
  RDSA & Gaussian &${Q}_2+D$\\
  RDSA & Uniformly Spherical &$\dfrac{p}{p+2}{Q}_2+D$\\
  \bottomrule %添加表格底部粗线
  \end{tabular}
  \label{table22}

\end{table}

 From the two tables, we notice that when dimension $p$ is sufficiently large, $\upphi, \upsilon$ of uniformly spherical are identical to those of the Gaussian distribution. Also, as the order $d$ of $x^d$ in the U-shape distribution increases, $\upxi^2$ of that distribution is close to that of Bernoulli as well (the curve of the U-shape converges to two  single points $\pm1$). Next, we are going to choose Bernoulli and Gaussian, respectively to be the representative distribution of SPSA and RDSA and make comparisons between them.

\textbf{Proposition 3}: If $2\bm{u}^T_1\bm{Su}_1+\bm{u}^T_1\bm{Su}_2\geq0$ holds, then SPSA with Bernoulli distributed perturbation has a smaller MSE value than RDSA with Gaussian perturbation. 

\textbf{Corollary 3A:} If $|2\bm{u}_1|>|\bm{u}_2|$ holds component-wise and $S_{ij}\geq0$, then Proposition 3 holds.

\noindent {\it{Proof:}} $\bm{S}=(a\bm{H}(\bm{\uptheta}^*)-\upbeta_{+}\bm{I}/2)^{-2}$ is a positive-definite matrix by the fact that  $\bm{H}(\bm{\uptheta}^*)$ is diagonalizable and $\upbeta_{+}<2\text{min}\uplambda_i(a\bm{H}(\uptheta^*))$. By matrix analysis, for any $\bm{u}_1,\bm{u}_2\in\mathbb{R}^p$:
$$|\bm{u}^T_1\bm{S}\bm{u}_2|\leq\sqrt{(\bm{u}^T_1\bm{Su}_1)(\bm{u}^T_2\bm{Su}_2)},$$

So we would have, 
\begin{align}
\begin{split}
    &2\bm{u}^T_1\bm{Su}_1+\bm{u}^T_1\bm{Su}_2\\
    &\geq 2\bm{u}^T_1\bm{Su}_1-\sqrt{(\bm{u}^T_1\bm{Su}_1)(\bm{u}^T_2\bm{Su}_2)}.
\label{eq:so}
\end{split}
\end{align}

If $|2\bm{u}_1|\geq|\bm{u}_2|$ component-wise and $S_{ij}\geq0$, then the r.h.s of (\ref{eq:so}) would be larger than 0, so that $2\bm{u}^T_1\bm{Su}_1+\bm{u}^T_1\bm{Su}_2\geq0$, leading to,
\begin{align*}
    (3\bm{u}_1+\bm{u}_2)^T\bm{S}(3\bm{u}_1+\bm{u}_2)\geq(\bm{u}_1+\bm{u}_2)^T\bm{S}(\bm{u}_1+\bm{u}_2).
\end{align*}

Consequently, the MSE of SPSA with Bernoulli-distributed perturbation will be lower than that of RDSA with Gaussian-distributed perturbation.

\textbf{Corollary 3B:} If the loss function does not have cross-third derivative term, i.e. $\bm{u}_2=\bm{0}$ along with $\upxi^2<\upphi, \uprho<1$, then SPSA has a smaller MSE value than RDSA for any valid perturbation.

\noindent {\it{Proof:}}
When $\bm{u}_2=\bm{0}$, we would have:
\begin{align}
\begin{split}
     &\text{MSE}_{\text{SPSA}}=\upxi^4\bm{u}^T_1\bm{Su}_1+\uprho^2D,\\ &\text{MSE}_{\text{RDSA}}=\upphi^2\bm{u}^T_1\bm{Su}_1+D,
\end{split}
    \label{eq:zhiqian}
\end{align}

Then under the condition that $\upxi^2<\upphi, \uprho<1$, SPSA would have a smaller MSE value.

\subsubsection{Application on the skewed-quartic function}

Consider the skewed quartic loss function:
$$L(\bm{\uptheta})=\bm{\uptheta}^T\bm{B}^T\bm{B\uptheta}+0.1\sum\limits^p_{i=1}(\bm{B\uptheta})^3_i+0.01\sum\limits^p_{i=1}(\bm{B\uptheta})^4_i.$$
where $\bm{\uptheta}=(t_1,t_2,...,t_p)^T$. After calculation we obtain that:

$\bm{u}_1=ac^2(L^{'''}_{111}(\bm{\uptheta}),...,L^{'''}_{ppp}(\bm{\uptheta}))^T/6|_{\bm{\uptheta}^*=\bm{0}}$=
\begin{small}
\begin{align*}
    \dfrac{ac^2}{6}\left(                
  \begin{array}{l}  
    \dfrac{0.6}{p^3}+\dfrac{0.01}{p^4}\times 24\sum\limits^p_{i=1}t_i\\  
    \dfrac{0.6}{p^3}\times2+\dfrac{0.01}{p^4}\times 24\times(\sum\limits^p_{i=1}t_i+\sum\limits^p_{i=2}t_i)\\ 
    \dfrac{0.6}{p^3}\times3+\dfrac{0.01}{p^4}\times 24\times(\sum\limits^p_{i=1}t_i+\sum\limits^p_{i=2}t_i+\sum\limits^p_{i=3}t_i)\\
    ...\\
    \dfrac{0.6}{p^3}\times p+\dfrac{0.01}{p^4}\times 24\times(\sum\limits^p_{i=1}t_i+\sum\limits^p_{i=2}t_i+\sum\limits^p_{i=3}t_i+...+t_p)
  \end{array}
\right)_{\bm{\uptheta}^*=\bm{0}}
=
\dfrac{0.1\times ac^2}{p^3} \left(                
  \begin{array}{c}  
   1\\
   2\\
   3\\
   ...\\
   p-1\\
   p
   \end{array}
\right)
\end{align*}
\end{small}

$\bm{u}_2=ac^2(3\sum\limits^p_{j\neq1}L^{'''}_{jj1}(\bm{\uptheta}),...,3\sum\limits^p_{j\neq p}L^{'''}_{jjp}(\bm{\uptheta}))^T/6|_{\bm{\uptheta}^*=\bm{0}}$=

\begin{align*}
%\begin{footnotesize}
    \dfrac{ac^2}{2}\left(                
  \begin{array}{l}  
   (p-1)\left(\dfrac{0.6}{p^3}+\dfrac{0.01}{p^4}\times24\sum\limits^p_{i=1}t_i\right)\\
   (p-1)\left(\dfrac{0.6}{p^3}+\dfrac{0.01}{p^4}\times24\sum\limits^p_{i=1}t_i\right)+(p-2)\left(\dfrac{0.6}{p^3}+\dfrac{0.01}{p^4}\times24\sum\limits^p_{i=2}t_i\right)\\
   ...\\
   (p-1)\left(\dfrac{0.6}{p^3}+\dfrac{0.01}{p^4}\times24\sum\limits^p_{i=1}t_i\right)+...+\dfrac{0.6}{p^3}+\dfrac{0.01}{p^4}\times24\sum\limits^p_{i=p-1}t_i
   \end{array}
\right)
%\end{footnotesize}
\end{align*}

\begin{align*}
   =\dfrac{0.3\times ac^2}{p^3}\left(                
  \begin{array}{c}  
   p-1\\
   (p-1)+(p-2)\\
   ...\\
   (p-1)+(p-2)+...+2+1
   \end{array}
\right)
\end{align*}

$\bm{H}(\bm{\uptheta}^*)=(L^{''}_{ij}(\bm{\uptheta}))_{ij}|_{\bm{\uptheta}^*=\bm{0}}$=
\begin{align*}
    \left(                
  \begin{array}{ccccc}  
   \dfrac{2}{p^2} & \dfrac{2}{p^2}& \dfrac{2}{p^2} &...&\dfrac{2}{p^2}\\
   \dfrac{2}{p^2}& \dfrac{2}{p^2}\times2 &\dfrac{2}{p^2}\times2 &...&\dfrac{2}{p^2}\times2\\
   \dfrac{2}{p^2}& \dfrac{2}{p^2}\times2 &\dfrac{2}{p^2}\times3&...&\dfrac{2}{p^2}\times3\\
   \dfrac{2}{p^2}\times& \dfrac{2}{p^2}\times2 &\dfrac{2}{p^2}\times3&...&\dfrac{2}{p^2}\times4\\
   ... \\
   \dfrac{2}{p^2}& \dfrac{2}{p^2}\times 2 &\dfrac{2}{p^2}\times3&...&\dfrac{2}{p^2}\times p
   \end{array}
\right)=
   \dfrac{2}{p^2} \left(                
  \begin{array}{cccccc}  
   1 &1 &1 &1&...&1\\
   1 &2 &2 &2&... &2\\
   1 &2 &3 &3&...&3\\
   1 &2 &3 &4&...&4\\
   ...\\
   1&2&3&4&...&p
   \end{array}
\right)
\end{align*}

In many practical cases, we set $\upalpha$ in the gain sequence $a_k=a/(k+1)^{\upalpha}$ less than 1, so the $\upbeta^{+}=0$. Then

$\bm{S}=(a\bm{H}(\bm{\uptheta}^*)-\dfrac{\upbeta^{+}}{2}\bm{I})^{-2}$=
\begin{align*}
   \dfrac{p^4}{4a^2} \left(                
  \begin{array}{ccccccccc}  
   5 &-4 &1 &0&0 &...&0&0&0\\
   -4 &6 &-4 &1&0&... &0&0&0\\
   1 &-4 &6 &-4&1&...&0&0&0\\
   0 &1 &-4 &6&-4&...&0&0&0\\
   ...\\
   0&0&0&0&...&1&-4&6&-3\\
   0&0&0&0&...&0&1&-3&2
   \end{array}
\right)
\end{align*}

Then we could calculate the key terms:
\begin{align*}
    \begin{split}
        &\bm{u}^T_1\bm{Su}_1=\dfrac{0.01\times c^4}{4a^2p^2};\\
        &\bm{u}^T_1\bm{Su}_2=0;\\
        &\bm{u}^T_2\bm{Su}_2=\dfrac{0.09\times c^4(p-1)}{4a^2p^2}.
    \end{split}
\end{align*}
So, we will have:
\begin{align}
    2\bm{u}^T_1\bm{Su}_1+\bm{u}^T_1\bm{Su}_2=\dfrac{0.01\times c^4}{2a^2p^2}.
    \label{eq:ai}
\end{align}

Since $p$ represents an integer larger than 1, the above term is always positive. According to Proposition 3, it is known that the asymptotic MSE value of SPSA with Bernoulli is smaller than the MSE for RDSA with Gaussian (or other perturbation with fourth moment of its distribution being 3 and the product of two separate second moments being 1). The corresponding numerical experiment is shown in the next section.

Moreover, we can obtain the MSE value for each case:
\begin{align*}
    \begin{split}
        &Q_1=(\bm{u}_1+\bm{u}_2)^T\bm{S}(\bm{u}_1+\bm{u}_2)=\dfrac{0.09c^4p-0.08c^4}{4a^2p^2},\\
        &Q_2=(3\bm{u}_1+\bm{u}_2)^T\bm{S}(3\bm{u}_1+\bm{u}_2)=\dfrac{0.09c^4}{4a^2p},\\
        &D=\dfrac{\sigma^2}{8c^2}\sum\limits_{i}\dfrac{1}{\uplambda_i}\;\;\;\; (\uplambda_i\;\text{are\;the\;eigenvalues\;of\;matrix\;}a\bm{H}(\bm{\uptheta}^*)).
    \end{split}
\end{align*}

From Table 2, we find that the asymptotic MSE value of Bernoulli will always be smaller than that of Gaussian since $Q_1<Q_2$. Also, we find that the MSE value of Uniform Spherical will be smaller than Gaussian with the coefficient in front of $Q_2$ being less than 1. For the U-shape, there is no certain conclusion for the relative MSE value since one contribution $D$ is going down while the other one $Q_1$ is going up.

\subsection {A More Straightforward Look}
We can have a more straightforward understanding of the comparison of $Q_1$ and $Q_2$ by looking in a 2-dim way. Let $S, u_1$ and $u_2$ be scalars, set $x=\sqrt{S}u_1, y=\sqrt{S}u_2$, then consider $(3x+y)^2-(x+y)^2$:
\begin{align}
    &(3x+y)^2-(x+y)^2\\
    &=9x^2+6xy+y^2-x^2-2xy-y^2=8x^2+4xy.
\label{eq:haha}
\end{align}

Under the following conditions, we will have (14) larger than 0. 
$$
\left\{
\begin{array}{lr}
y \geq -2x\;\;\;\text{if\;} x>0\\
y< -2x\;\;\;\text{if\;} x<0
\end{array}
\right.
$$

The shadow region in the following figure shows the area satisfying the above condition, which indicates a large probability, around 5/6 (more accurately,  1-2*arctan (0.5)) of $Q_2-Q_1>0$. In other words, there is a probability around 5/6 that SPSA with Bernoulli outperforms RDSA with Gaussian.

\begin{figure}[H]
		\centering
		\includegraphics[width=10cm,height=10cm]{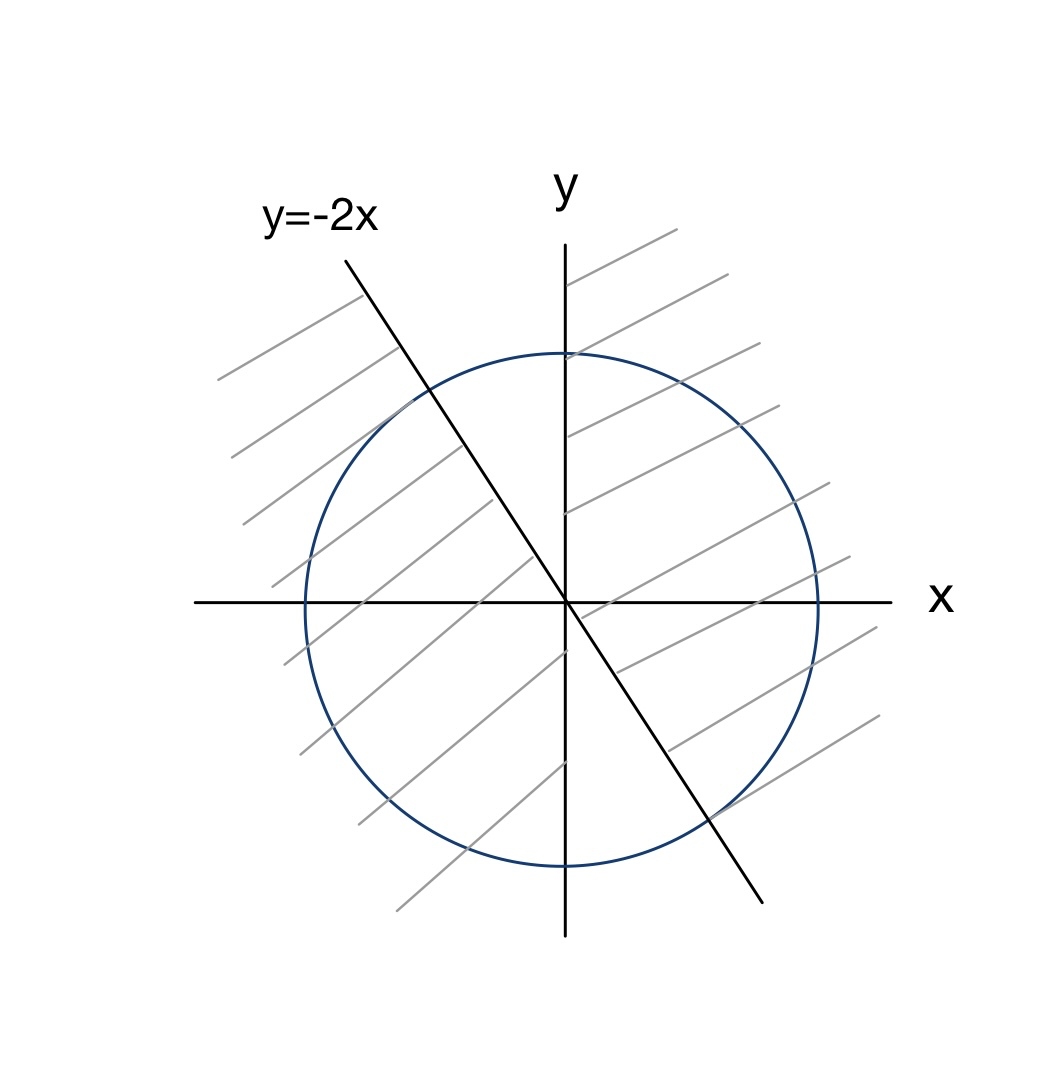}
		\caption{Region of equation (14) being larger than 0 }
	\end{figure}

Then we are going to consider the high-dim case. Denote:
\begin{align*}
    \begin{split}
        &{Q}_1=(\bm{u}_1+\bm{u}_2)^T\bm{S}(\bm{u}_1+\bm{u}_2),\\ &{Q}_2=(3\bm{u}_1+\bm{u}_2)^T\bm{S}(3\bm{u}_1+\bm{u}_2).
    \end{split}
\end{align*}

We would like to compare $Q_2, Q_1$. Similar to the scalar case, set $\bm{x}=\bm{S}^{1/2}\bm{u}_1, \bm{y}=\bm{S}^{1/2}\bm{u}_2.$ Then consider:
\begin{align}
    &\bm{z}=Q_2-Q_1\\
    &=(3\bm{x}+\bm{y})^T(3\bm{x}+\bm{y})-(\bm{x}+\bm{y})^T(\bm{x}+\bm{y})\\
    &=4(2\bm{x}^T\bm{x}+\bm{x}^T\bm{y})\\
    &=4\bm{x}^T(2\bm{x}+\bm{y})
    \label{eq:haha2}
\end{align}

Here we consider a simple case to have a straightforward understanding. Set $\bm{x},\bm{y}$ to be two independent random variables with uniform symmetrical distribution (i.e. $\bm{x}=(x_1,x_2,...,x_p)^T, \bm{y}=(y_1,y_2,...,y_p)^T,\; x_i\stackrel{\text{i.i.d}}\sim U(-a,a),  y_i\stackrel{\text{i.i.d}}\sim U(-b,b)$) and $\bm{S}$ to be the identical matrix. Then we have:
\begin{align*}
    &E(\bm{x}^T\bm{x})=E\left(\sum\limits^p_{i=1} x^2_i\right)=\sum\limits^p_{i=1} E(x^2_i)=\sum\limits^p_{i=1} Var(x_i)=\dfrac{a^2p}{3}.\\
    &E(\bm{x}^T\bm{y})=E\left(\sum\limits^p_{i=1} x_iy_i\right)=\sum\limits^p_{i=1} E(x_i)E(y_i)=0.\;\;\;(\text{due\;to\;the\;independence})\\
    &E(\bm{z})=8E(\bm{x}^T\bm{x})+4E(\bm{x}^T\bm{y})=\dfrac{8a^2p}{3}.
\end{align*}

Furthermore, we use one-sided version of Chebyshev inequality to see the probability that $\bm{z}\leq0$, which indicates $Q_2\leq Q_1$. The inequality states that $P(X\geq c)\leq Var(X)/[Var(X)+c^2]$. In our case, the inequality is in the form of: 
\begin{align}
    P[(-\bm{z}+E(\bm{z}))\geq E(\bm{z})]&\leq \dfrac{Var[-\bm{z}+E(\bm{z})]}{Var[-\bm{z}+E(\bm{z})]+[E(\bm{z})]^2}\\
    &= \dfrac{Var(\bm{z})}{Var(\bm{z})+[E(\bm{z})]^2}
    \label{eq:kun}.
\end{align}

To simplify the calculation, we set $x_i\stackrel{\text{i.i.d}}\sim U(-a,a), y_i\stackrel{\text{i.i.d}}\sim U(-a,a)$. Then we have:
\begin{align*}
    Var(\bm{z})&=Var(8\bm{x}^T\bm{x})+Var(4\bm{x}^T\bm{y})+Cov(8\bm{x}^T\bm{x},4\bm{x}^T\bm{y})\\
    &=64Var\left(\sum\limits^p_{i=1}x^2_i\right)+16Var\left(\sum\limits^p_{i=1}x_iy_i\right)=\dfrac{64pa^4}{5}+\dfrac{16pa^4}{9}=16\times\dfrac{41pa^4}{45}.
\end{align*}
Plug it into the inequality (20), it looks like:
\begin{align}
    P(\bm{z}\leq 0)) &\leq \dfrac{Var(\bm{z})}{Var(\bm{z})+[E(\bm{z})]^2}\\
    & = \dfrac{16\times\dfrac{41pa^4}{45}}{16\times\dfrac{41pa^4}{45}+\dfrac{64p^2a^4}{9}}\\
    &= \dfrac{41}{41+20p}.
    \label{eq:ee}
\end{align}

Inequality (23) shows that as the dimension $p$ increases, the probability of $\bm{z}\leq 0$ decreases. Here we make a simulation by setting a range of dimension from 1 to 10 and then take 100,000 independent trials to calculate the probability, i.e. $P(\bm{z}\leq0)$ with $\bm{S}=\bm{I};\; \bm{u}_{1i}\stackrel{\text{i.i.d}}\sim U(-100,100);\; \bm{u}_{2i}\stackrel{\text{i.i.d}}\sim U(-100,100),\; i=1,2,...,p$. The result is shown in the following Table 3, demonstrating a decreasing trend when the dimension gets large. Also, it shows consistency with the conclusion in the scalar case, namely, when $p=1, P(z\leq0)\approx\dfrac{1}{6}$.

\begin{table}[ht!]
 \renewcommand\arraystretch{1.5}
 \centering
 \caption{Relationship Between Dimension and $P(\bm{z}\leq0)$}
  \begin{tabular}{cccccc}
  \toprule  %添加表格头部粗线
  \textbf{Dim} & 1&2&3&4&5\\
  \midrule  %添加表格中横线
  $\bm{P}(\bm{z}\leq0)$& 0.12546 & 0.0336& 0.00932& 0.00252& 0.00071\\
  \midrule  %添加表格中横线
  \textbf{Dim} & 6&7&8&9&10\\
  \midrule  %添加表格中横线
  $\bm{P}(\bm{z}\leq0)$& 0.00024 & $9\times 10^{-5}$ & $3\times 10^{-5}$ & 0& 0\\
  \bottomrule %添加表格底部粗线
  \end{tabular}
  \label{table3}
\end{table}

As Table 3 indicates, if we consider those high-dimension cases, there will be a tiny probability that $Q_2\leq Q_1$, i.e. SPSA with Bernoulli will highly probably beat RDSA with Gaussian from the perspective of a smaller MSE. However, it should be noticed that the above case is based on a simple simulated setting, which aims at providing an intuitive view at the comparison of Gaussian and Bernoulli. 

In the next section, we will show how Proposition 3 works more specifically by some numerical experiments.

\section{Numerical Experiment}
\subsection{Numerical Example 1}
Consider the loss function \cite{558808}: $$L(\uptheta)=||\uptheta||^2+\sum\limits_{i=1}^pe^{t_i/p},$$
where $\bm{\uptheta}=[t_1,t_2,...,t_p]^T$ and the minimum occurs at $\bm{\uptheta}^*$ with each component $-0.033$, such that $L(\bm{\uptheta}^*)=29.99994$. This is a convex function with third-derivative being nonzero and cross third-derivative being zero. Namely, in our case, we have $\bm{u}_2=\bm{0}$ and $\bm{u}_1>\bm{0}$. So according to the Corollary in Proposition 3, SPSA with Bernoulli has lower MSE than RDSA with Gaussian.

Let $p=30$ and take 100 independent trials, each with 3000 iterations. Results are presented in Fig 2 and Table 4. We assume the noise term has distribution $N(0,0.01)$ and gain sequence as $a_k=0.05/(k+1)^{0.602},\, c_k=0.3/(k+1)^{0.101}$. Set the initial point $\hat{\bm{\uptheta}}_0=[1,1,...,1]^T$.  The results are shown below. The lowest MSE is highlighted in bold in Table 4 and other numerical results below.
 \begin{figure}[H]
		\centering
		\includegraphics[width=12cm,height=8cm]{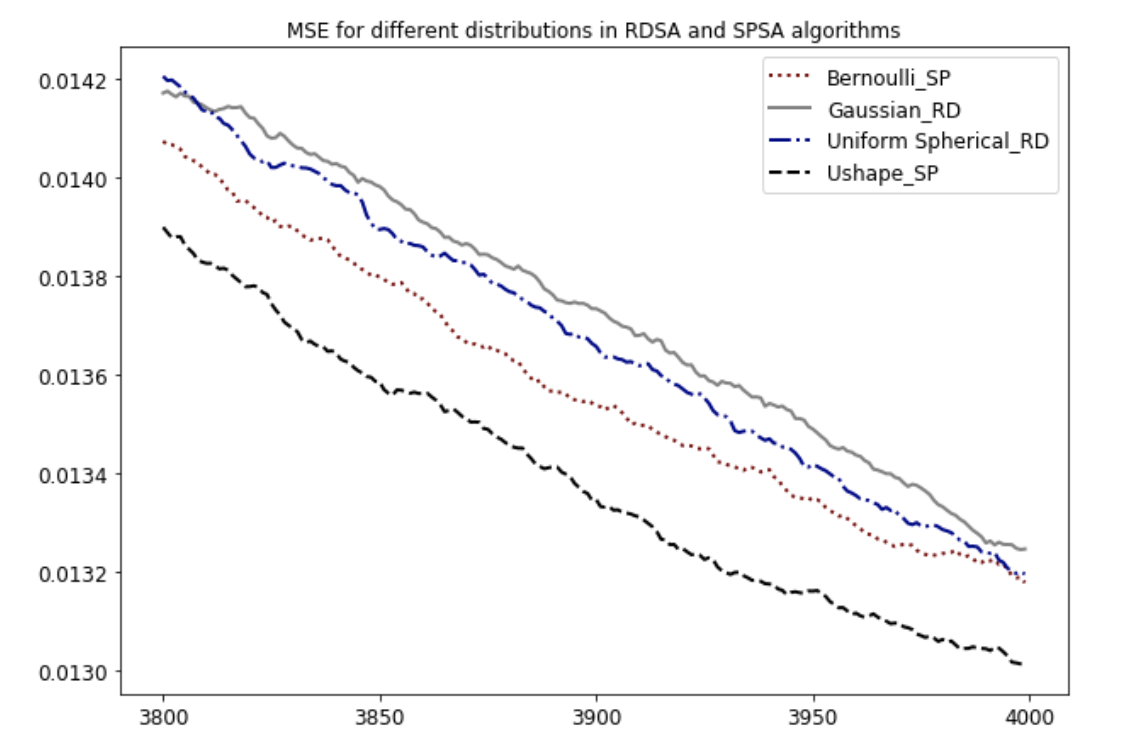}
		\caption{Mean Square Error, computed as an average over 100 trials, plotted the last 200 iterations for different random perturbation distributions}
	\end{figure}

\begin{table}[ht!]
 \renewcommand\arraystretch{1.5}
 \centering
 \caption{MSE and confidence interval in SPSA and RDSA with different perturbations with same gain sequences in all algorithms}
  \begin{tabular}{ccc}
  \toprule  %添加表格头部粗线
  Perturbation & MSE & 95\% CI\\
  \midrule  %添加表格中横线
  Bernoulli SP& 0.01318 & [0.01281, 0.01354]\\
  \textbf{U-shape}($\bm{x^{10}}$) \textbf{SP} &\textbf{0.01301} & \textbf{[0.01266, 0.01336]} \\
  Gaussian {\text{RD}} &0.01325 &    [0.01290, 0.01359]\\
  Uniformly Spherical {\text{RD}} &0.01319 & [0.01285, 0.01355]\\
  \bottomrule %添加表格底部粗线
  \end{tabular}
  \label{table4}
\end{table}

The experimental result shows consistency with the Proposition 3. The confidence intervals show overlap among the four and then we conduct the two-sample $t$-test between the Bernoulli and Gaussian cases and obtain the $p$-value 0.35982, which is larger than the common threshold, 0.05, indicating that the difference is not statistically significant. The reason is that in this case, the $D$ value (shown in the Table \ref{table22}) is much larger than the term $\bm{u}_1^T\bm{Su}_1$ so that it dominates the whole MSE value. So the MSE of SPSA with Bernoulli is slightly smaller than that of RDSA with Gaussian. Moreover, we notice that the MSE of U-Shape distributed perturbation is the smallest one which can also be explained by the 0.9 in front of the value $D$.

\subsection{Numerical Example 2}
Consider the Ackley Function with dimension $p=30$:

\begin{align*}
%\begin{split}
    L(\bm{\uptheta})&=-a\times\text{exp}\left(-b\sqrt{\dfrac{1}{30}\sum\limits^d_{i=1}t^2_i}\right)-\text{exp}\left(\dfrac{1}{30}\sum\limits^d_{i=1}\text{cos}(ct_i)\right)+a+e,
   % \end{split}
\end{align*}
where $a=20, b=0.2, c=2\pi,\,\bm{\uptheta}=[t_1,t_2,...,t_{30}]$ and $L(\bm{\uptheta})=0$ at $\bm{\uptheta}^*=[1,1,...,1]^T.$ Here, we assume the noise term has distribution $N(0,0.01)$ and gain sequence as $a_k=0.02/(k+1+A)^{0.602},\;A=10,\, c_k=0.2/(k+1)^{0.101}$. Set the initial point $\hat{\bm{\uptheta}}_0=[0.2,0.2,...,0.2]^T$. 

In this case, there are no cross-third derivatives, i.e. \,$\bm{u}_2=\bm{0}$. By Proposition 3, SPSA with Bernoulli would have smaller MSE value than RDSA with Gaussian. Take 100 independent trials, each with 5000 iterations. The results are shown in Fig 3 and Table 5:
 \begin{figure}[H]
		\centering
		\includegraphics[width=12cm,height=8cm]{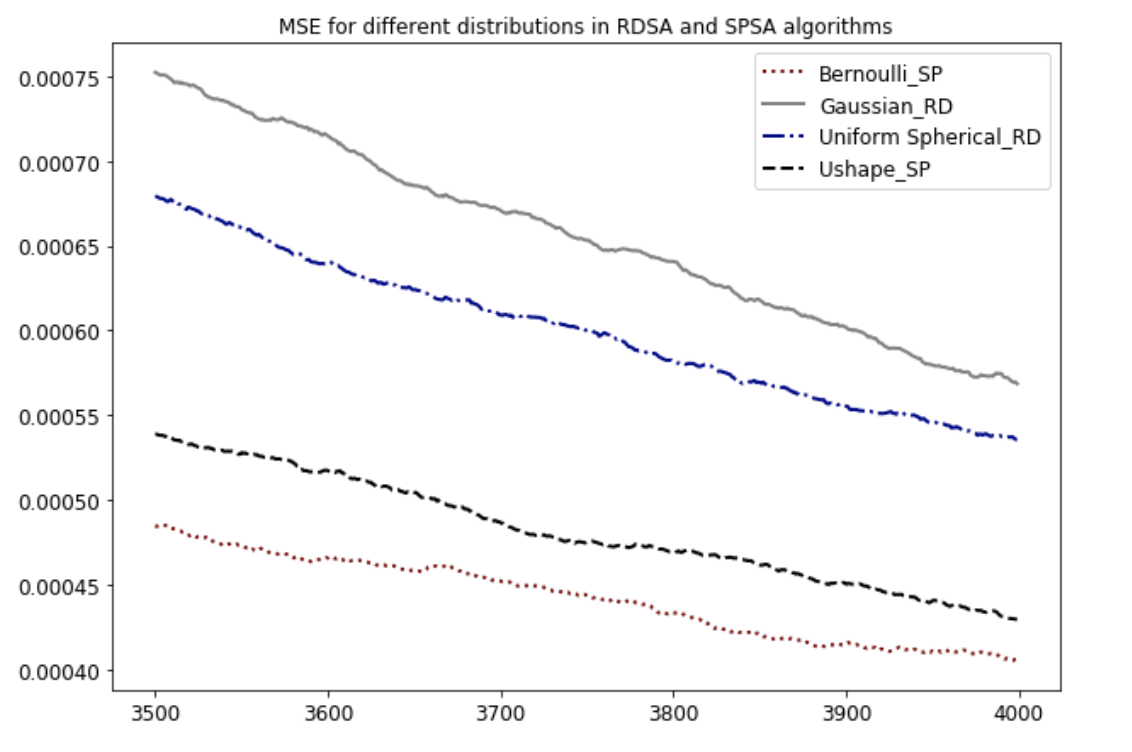}
		\caption{Mean Square Error, computed as an average over 100 trials, plotted the last 500 iterations for different random perturbation distributions}
	\end{figure}

\begin{table}[ht!]
 \renewcommand\arraystretch{1.5}
 \centering
 \caption{MSE and its Confidence Interval in SPSA and RDSA with different perturbations with same gain sequences in all algorithms}
  \begin{tabular}{ccc}
  \toprule  %添加表格头部粗线
  Perturbation & MSE & 95\% CI \\
  \midrule  %添加表格中横线
  \textbf{Bernoulli} \textbf{SP}& \textbf{0.00041} & \textbf{[0.00038, 0.00044]}\\
  U-shape ($x^{10}$) SP& 0.00043 &  [0.00041, 0.00045]\\
  Uniformly Spherical RD & 0.00054 & [0.00051, 0.00056]\\
  Gaussian RD & 0.00057 &   [0.00054, 0.00060]\\
  \bottomrule %添加表格底部粗线
  \end{tabular}
  \label{table5}
\end{table}

The MSE value of SPSA with Bernoulli case is lower than that of RDSA with Gaussian case, which is accordant with the Proposition 3. Notice that the confidence interval of the two cases does not overlap and we take the two-sample $t$-test between them and gain the $p$-value $5.96\times10^{-13}$, which indicates that the difference is significant. 

\subsection{Numerical Example 3}

To give a more specific example. Consider the skewed-quartic function: $$L({\bm{\uptheta}})={\bm{\uptheta}}^T{\bm{B}}^T{\bm{B\uptheta}}+0.1\sum\limits^p_{i=1}({\bm{B\uptheta}})^3_i+0.01\sum\limits^p_{i=1}({\bm{B\uptheta}})^4_i.$$

\subsubsection{Optimal Gain Sequence Searching}

We consider searching the `optimal gain sequence' for $a_k$ and $c_k$  in the gain sequence with the smallest MSE value, i.e. set $a_k=a/(k+A+1)^\alpha,\; c_k=c/(k+1)^\gamma$ and keep $A=10$. Then set the search range of $a$ and $c$ to be 0.1 to 1, with interval 0.02 and pick the sequence leading to the lowest MSE value. The `optimal gain sequence' for each case has been listed in Table 6. Particularly, in order to satisfy that $\bm{\upmu}\neq \bm{0}$ in (\ref{eq:mu}), we set $\upalpha=0.606, \upgamma=0.101$ in the gain sequence.

\begin{table}[ht!]
 \renewcommand\arraystretch{1.5}
 \centering
 \caption{Optimal Gain Sequence for Algorithms with distinct Perturbations}
  \begin{tabular}{ccc}
  \toprule  %添加表格头部粗线
  Perturbation & $a_k$ & $c_k$ \\
  \midrule  %添加表格中横线
  Bernoulli SP &$0.12/(k+11)^{0.606}$ & $0.8/(k+1)^{0.101}$\\
  U-shape ($x^{10}$) SP& $0.1/(k+11)^{0.606}$& $0.48/(k+1)^{0.101}$\\
  Uniformly Spherical RD & $0.1/(k+11)^{0.606}$ & $0.42/(k+1)^{0.101}$\\
  Gaussian RD & 0.1/$(k+11)^{0.606}$ & $0.58/(k+1)^{0.101}$\\
  \bottomrule %添加表格底部粗线
  \end{tabular}
  \label{table6}

\end{table}

\subsubsection{Further Comparison}

In order to conduct a further comparison of the performance of RDSA and SPSA algorithm, as well as work in connection with our previous theoretical calculation, we are going to 
apply the optimal gain sequence for each distribution to the other three distributions and see the result. That is, in 1)-4) below, we, in turn, pick one of the pairs {\{$a_k,c_k$\}} in Table 6 and run all four methods at the chosen gain sequence. This will demonstrate each method in its best configuration. Then we mainly focus on the Bernoulli and Gaussian to see their comparison.

\begin{itemize}
  \item [1)] 
  Apply the optimal gain sequence for Bernoulli, i.e. set $a_k=0.12/((k+11)^{0.606}), c_k=0.8/((k+1)^{0.101}$ for the four methods. The result is shown in Fig 4 and Table 7.
  \begin{figure}[H]
		\centering
		\includegraphics[width=13cm,height=9.5cm]{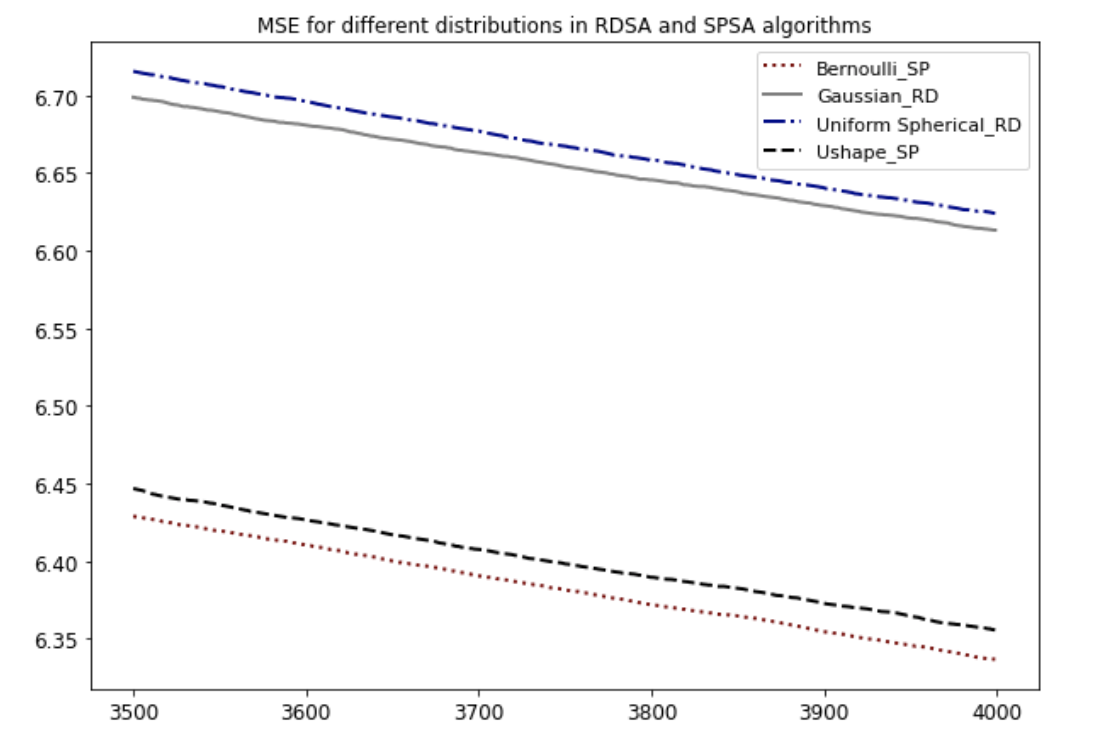}
		\caption{Mean Square Error, computed as an average over 100 trials each of 4000 iterations with the optimal gain for Bernoulli distributed perturbation in SPSA, plotted the last 500 iterations for different random direction distributions}
	\end{figure}
	\begin{table}[ht!]
 \renewcommand\arraystretch{1.5}
 \centering
 \caption{MSE for Skewed Quartic Loss Function and its Confidence Interval\\(Gain Sequence $a_k=0.12/((k+11)^{0.606}), c_k=0.8/((k+1)^{0.101}$)}
  \begin{tabular}{ccc}
  \toprule  %添加表格头部粗线
  Perturbation & MSE & 95\% CI \\
  \midrule  %添加表格中横线
  \textbf{Bernoulli} \textbf{SP}& \textbf{6.3368} &  \textbf{[6.0959, 6.5778]}\\
  U-shape ($x^{10}$) SP& 6.3557 &  [6.0536, 6.6579]\\
  Uniformly Spherical RD & 6.6240 & [6.3064, 6.9416]\\
  Gaussian RD & 6.6132 &   [6.2637, 6.9626]\\
  \bottomrule %添加表格底部粗线
  \end{tabular}
  \label{table7}
\end{table}

 We find that the Bernoulli performs the best, then the U-shape, Uniformly Spherical and Gaussian being the worst one. It is accordant with the previous Table \ref{table22}. Also, we find that the 95\% confidence interval of Bernoulli and Gaussian does not overlap. The two-sample test between Bernoulli and Gaussian produces the $p$-value 0.0040573, indicating that the difference is significant.
 
  \item [2)]
   Apply the optimal gain sequence for U-shape, i.e. set $a_k=0.1/((k+11)^{0.606}), c_k=0.48/((k+1)^{0.101}$ for the four methods. The result is shown in Fig 5 and Table 8.
   \begin{figure}[H]
		\centering
		\includegraphics[width=12cm,height=9cm]{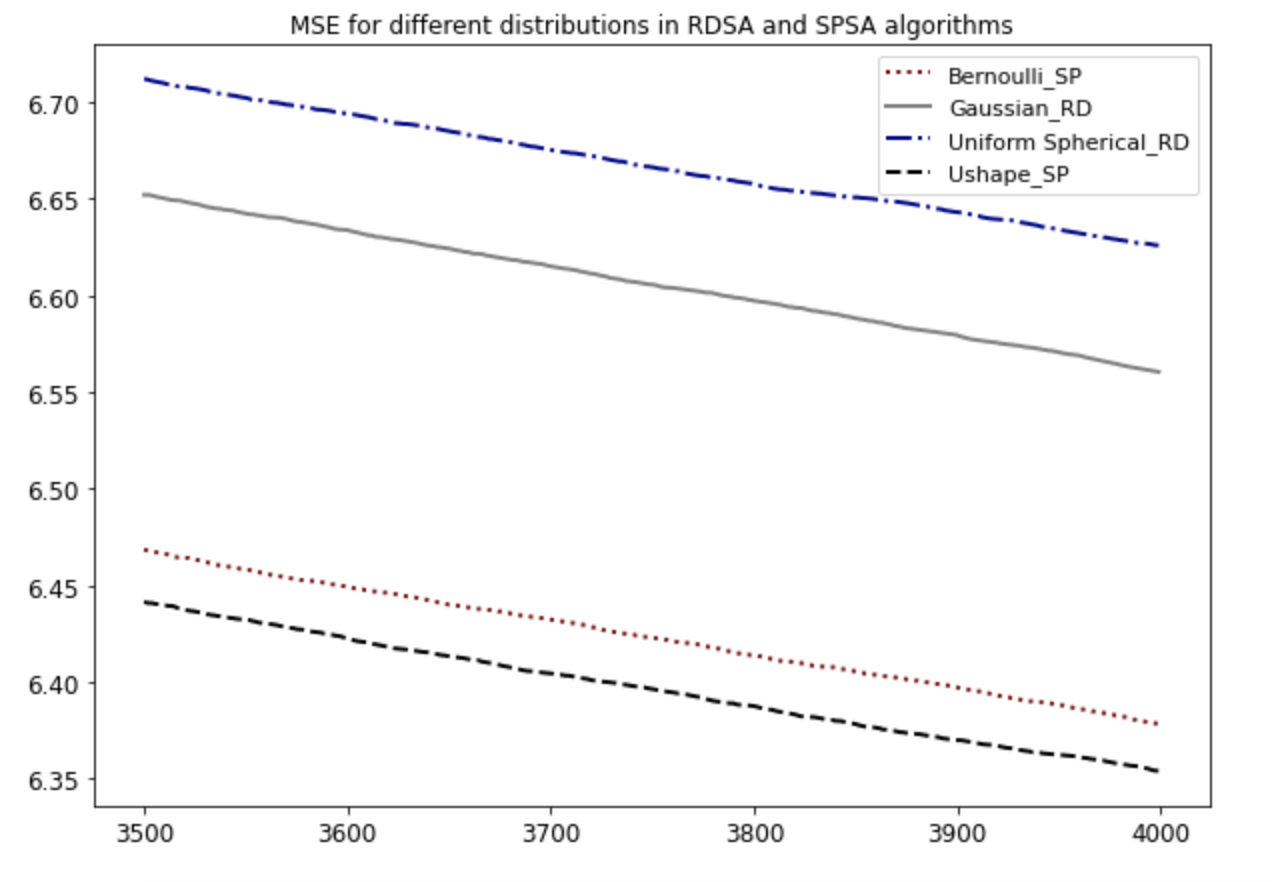}
		\caption{Mean Square Error, computed as an average over 100 trials each of 4000 iterations with the optimal gain for U-shape in SPSA, plotted the last 500 iterations for different random direction distributions}
	\end{figure}
	
		\begin{table}[ht!]
 \renewcommand\arraystretch{1.5}
 \centering
 \caption{MSE for Skewed Quartic Loss Function and its Confidence Interval\\ (Gain Sequence $a_k=0.1/((k+11)^{0.606}), c_k=0.48/((k+1)^{0.101}$)}
  \begin{tabular}{ccc}
  \toprule  %添加表格头部粗线
  Perturbation & MSE & 95\% CI \\
  \midrule  %添加表格中横线
  Bernoulli SP & 6.3784 &  [6.1038, 6.6531]\\
  \textbf{U-shape}($\bm{x^{10}}$) \textbf{SP}& \textbf{6.3540} &  \textbf{[6.0761, 6.6320]}\\
  Uniformly Spherical RD &6.6257 & [6.3150, 6.9365]\\
  Gaussian RD & 6.5604 &   [6.2803, 6.8406]\\
  \bottomrule %添加表格底部粗线
  \end{tabular}
  \label{table8}
\end{table}
   
  \item [3)]
  Apply the optimal gain sequence for Uniform Spherical, i.e. set $a_k=0.1/((k+11)^{0.606}), c_k=0.42 /((k+1)^{0.101}$ for the four methods. The result is shown in Fig 6 and Table 9.
  \begin{figure}[H]
		\centering
		\includegraphics[width=13cm,height=9cm]{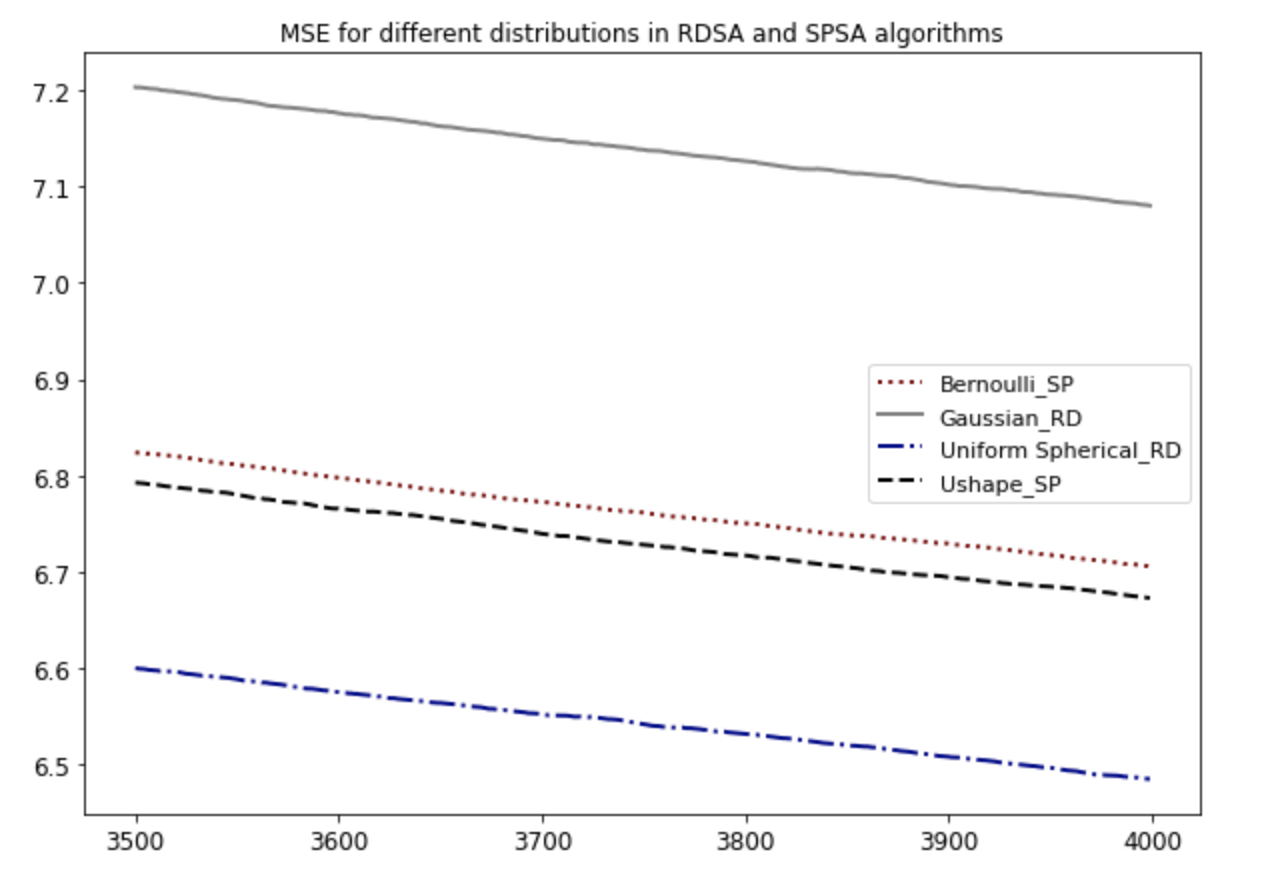}
		\caption{Mean Square Error, computed as an average over 100 trials each of 4000 iterations with the optimal gain for Uniform Spherical in RDSA, plotted the last 500 iterations for different random direction distributions}
	\end{figure}
		\begin{table}[ht!]
 \renewcommand\arraystretch{1.5}
 \centering
 \caption{MSE for Skewed Quartic Loss Function and its Confidence Interval\\(Gain Sequence $a_k=0.1/((k+11)^{0.606}), c_k=0.42 /((k+1)^{0.101}$)}
  \begin{tabular}{ccc}
  \toprule  %添加表格头部粗线
  Perturbation & MSE & 95\% CI \\
  \midrule  %添加表格中横线
  Bernoulli SP& 6.7323 &  [6.4523,7.0123]\\
  U-shape ($x^{10}$) SP& 6.7012 &  [6.4212,6.9812]\\
  \textbf{Uniformly Spherical} \textbf{RD} & \textbf{6.4436} & \textbf{[6.1636,6.7236]}\\
  Gaussian RD & 7.0789 &   [6.7889,7.3689]\\
  \bottomrule %添加表格底部粗线
  \end{tabular}
  \label{table9}
\end{table}
	
 We find that under such gain sequence, the Uniform Spherical performs the best and Bernoulli still outperforms Gaussian, consistent with the Table \ref{table22}. We then conduct the two-sample $t$-test between Bernoulli and Gaussian and obtain the $p$-value 0.121487, indicating that the difference is not statistically significant.
	
  \item [4)]
  Apply the optimal gain sequence for Gaussian, i.e. set $a_k=0.1/((k+11)^{0.606}), c_k=0.58 /((k+1)^{0.101}$ for the four methods. The result is shown in Fig 7 and Table 10.
   \begin{figure}[H]
		\centering
		\includegraphics[width=12.5cm,height=9.5cm]{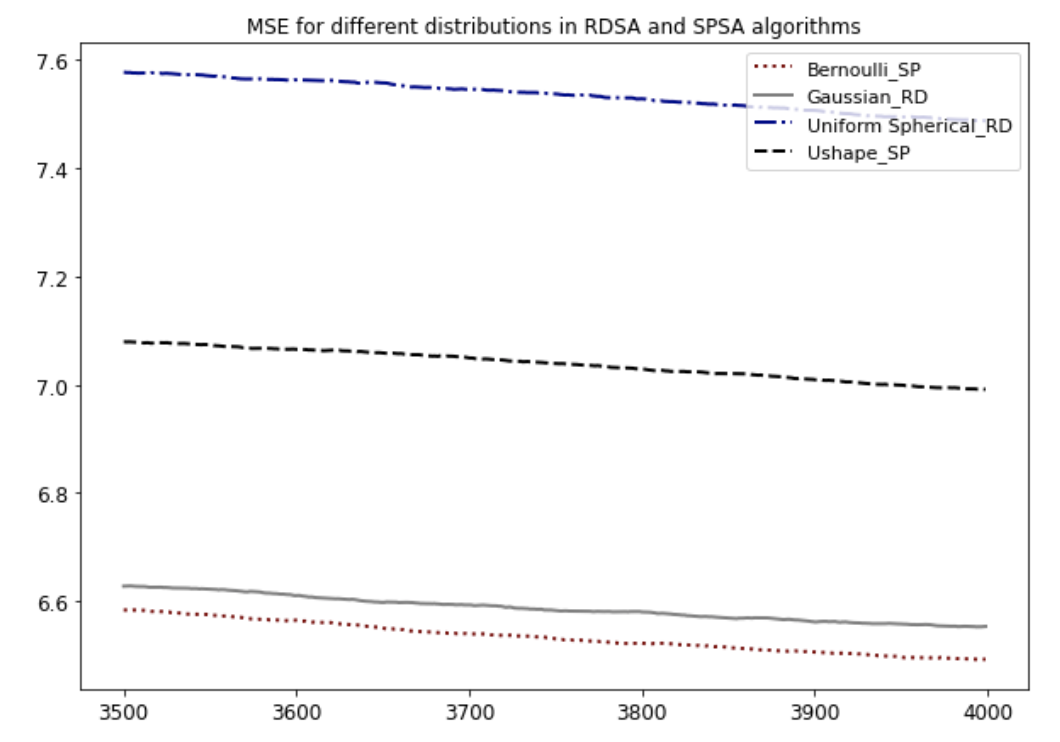}
		\caption{Mean Square Error, computed as an average over 100 trials each of 4000 iterations with the optimal gain for Gaussian distributed perturbation in RDSA, plotted the last 500 iterations for different random direction distributions}
	\end{figure}
	\begin{table}[ht!]
 \renewcommand\arraystretch{1.5}
 \centering
 \caption{MSE for Skewed Quartic Loss Function and its Confidence Interval\\(Gain Sequence $a_k=0.1/((k+11)^{0.606}), c_k=0.58 /((k+1)^{0.101}$)}
  \begin{tabular}{ccc}
  \toprule  %添加表格头部粗线
  Perturbation & MSE & 95\% CI \\
  \midrule  %添加表格中横线
  \textbf{Bernoulli} \textbf{SP}& \textbf{6.4928} &  \textbf{[6.2314, 6.7542]}\\
  U-shape ($x^{10}$) SP& 6.9915 &  [6.7374, 7.2456]\\
  Uniformly Spherical RD & 7.4873 & [7.1479, 7.8267]\\
  Gaussian RD & {6.5534} &  {[6.2897, 6.8171]}\\
  \bottomrule %添加表格底部粗线
  \end{tabular}
  \label{table10}
\end{table}

\end{itemize}

In this case, although the optimal gain sequence of Gaussian has been applied and the Gaussian reaches its lowest MSE, it is still smaller than the MSE of Bernoulli, which  can be a strong support for the theoretical conclusion shown in Proposition 3. The two-sample test between Bernoulli and Gaussian produces the $p$-value 0.2431523, indicating that the difference is not statistically significant.

From above four experiments, all the four algorithms show their best performance from the perspective of a lowest MSE when their optimal gain sequence is applied. The tables show that no matter in which case, SPSA with Bernoulli beats RDSA with Gaussian, showing consistency with Proposition 3.

\section{Conclusion}
This paper provides a strict derivation for several essential terms in the RDSA algorithm and presents a comparison of the asymptotic MSE between RDSA and SPSA. Furthermore, this paper presents specific conditions under which SPSA with Bernoulli perturbation works better in the MSE sense than RDSA with Gaussian perturbation, which may provide practitioners a general guidance in the practical cases. In particular, it is found that SPSA outperforms RDSA across the large majority of implementation cases with different gain sequences and types of loss functions.

There are several directions for future research. First, we can consider extra types of random directions with different distributions. Second, more work related to robustness and relative efficiency would be useful. Third, since this paper mainly focuses on Bernoulli and Gaussian cases, more general cases could be considered for the comparison of MSE between SPSA and RDSA. It would also be of interest to carry out a comparison such as that above for both second-order or non-continuous versions of SPSA and RDSA, such as presented in \cite{4982684}, \cite{8430974}, \cite{9029707}, \cite{Zhu2020EfficientIO}, and \cite{7544612}, and constant-gain (step size) versions, such as in \cite{KAJIWARA2018396}, \cite{https://doi.org/10.1002/rnc.5151}, and \cite{9092353}.

\bibliography{myref.bib}
\end{document}